\documentclass{ifacconf}

\makeatletter
\let\old@ssect\@ssect % Store how ifacconf defines \@ssect
\makeatother

\usepackage{natbib}
\usepackage{hyperref}

\makeatletter
\def\@ssect#1#2#3#4#5#6{%
  \NR@gettitle{#6}% Insert key \nameref title grab
  \old@ssect{#1}{#2}{#3}{#4}{#5}{#6}% Restore ifacconf's \@ssect
}
\makeatother

\usepackage{graphicx}      % include this line if your document contains figures
\usepackage{url}
\usepackage{enumerate}
\usepackage{commath}
\usepackage{graphicx}
\usepackage{amsmath}
\usepackage{breqn}
\usepackage{subfigure}
\usepackage{booktabs}
\usepackage{amsfonts}
\usepackage{amssymb}
\usepackage{xcolor}
\usepackage{mathrsfs}
\usepackage{accents}
\usepackage{thmtools}
\usepackage{thm-restate}
\usepackage{xcolor}
\usepackage{algpseudocode}
\usepackage{algorithm}
\usepackage{pgf}
\begin{document}
\begin{frontmatter}

\title{Domain-aware Control-oriented Neural Models for Autonomous Underwater Vehicles\thanksref{footnoteinfo}} 

\thanks[footnoteinfo]{The research described in this work is part of the Mathematics for Artificial Reasoning in Science Initiative at Pacific Northwest National Laboratory. It was conducted under the Laboratory Directed Research and Development Program at PNNL, a multiprogram national laboratory operated by Battelle for the U.S. Department of Energy under contract DE-AC05-76RL01830.}

\author[First]{Wenceslao Shaw Cortez} 
\author[First]{Soumya Vasisht} 
\author[First]{Aaron Tuor}
\author[First]{J{\'a}n Drgo{\v n}a}
\author[First]{Draguna Vrabie}

\address[First]{Pacific Northwest National Laboratory, 
   Richland, WA 99354 USA (e-mail: \{w.shawcortez, soumya.vasisht, aaron.tuor, jan.drgona, draguna.vrabie\}@pnnl.gov).}

\begin{abstract}                % Abstract of not more than 250 words.
Conventional physics-based modeling is a time-consuming bottleneck in control design for complex nonlinear systems like autonomous underwater vehicles (AUVs). In contrast, purely data-driven models, though convenient and quick to obtain, require a large number of observations and lack operational guarantees for safety-critical systems. Data-driven models leveraging available partially characterized dynamics have potential to provide reliable systems models in a typical data-limited scenario for high value complex systems, thereby avoiding months of expensive expert modeling time. In this work we explore this middle-ground between expert-modeled and pure data-driven modeling. We present control-oriented parametric models with varying levels of domain-awareness that exploit known system structure and prior physics knowledge to create constrained deep neural dynamical system models. We employ universal differential equations to construct data-driven blackbox and graybox representations of the AUV dynamics. In addition, we explore a hybrid formulation that explicitly models the residual error related to imperfect graybox models. We compare the prediction performance of the learned models for different distributions of initial conditions and control inputs to assess their accuracy, generalization, and suitability for control.
\end{abstract}


\end{frontmatter}
%===============================================================================

\section{Introduction}
Underwater vehicle control is a vital aspect of maritime applications, such as marine environmental monitoring, seabed mapping, geohazard assessment, and remote exploration and sensing \citep{wynn2014AUVappl, yu2002AUVappl}. Control design for autonomous underwater vehicles (AUVs) is a complex task due to the underactuated nonlinear dynamics and the extreme operating environment. Most applications are often restricted to simplistic control strategies that are lacking in speed, versatility, and reliable autonomy in extreme environments and situations.
Advanced control methods such as predictive or adaptive control require reliable system models for prediction and to evaluate system properties such as stability and safety. Classical system identification techniques that rely on transfer function learning or parameter estimation of linear system models~\citep{LJUNG20101,VANOVERSCHEE199475} are inefficient and non-trivial for complex systems. %This may drive more efforts towards ``model-free" approaches~\citep{}, but the successes of model-free approaches are limited for real world systems, with fewer guarantees on convergence time and poor data efficiency. \sv{cite}
% \textcolor{red}{Not sure if this is correct in general. Model-free methods like PPC/funnel-based control ensure safety and convergence. The issues with PPC are sensitivity to noise and ensuring input constraints}. 
% Model-based approaches have seen considerable success in making long-term decisions, mostly amplified by the recent advances in machine learning techniques and tools \textcolor{red}{not sure if we can justify this? Maybe remove this sentence?}. 

Methods that discover governing equations underlying a dynamical system directly from system measurement data \citep{brunton2016sindy, raissi2017MFinf} have revolutionized the domain of scientific machine learning and data-driven dynamics discovery. Physics-informed neural networks (PINNs) \citep{raissi2019PINN} incorporate physical domain knowledge in the form of soft constraints to solve or obtain parameters of nonlinear partial differential equations (PDEs) from observed noisy data using automatic differentiation for training the deep learning models. However, PINNs have been shown to fail for complex problems, making the loss function harder to optimize \citep{krishnapriyan2021PINNs}. 
Moreover, PINNs have been predominantly used for learning the solutions to PDEs for a single initial condition, and thus are not directly suitable for use in control that requires a generative model trained over a distribution of initial conditions.

Deep learning-based alternatives to traditional system identification and control have become a leading research avenue due to their rich expressive power and their ability to capture complex dynamical interactions from data.
Neural state-space models (NSSMs) represent an example of such methods that use discrete-time formulations~\citep{NIPS2018_8004}.  \cite{skomski2021NSSM} showed that NSSMs can embed domain knowledge, including structural priors and block nonlinearities. These models have been shown to work well in predictive control and differentiable programming settings \citep{drgona2021DPC}.
A continuous-time and depth alternative to popular deep learning models such as residual networks (ResNets) was introduced in \cite{chen2018NODE} called neural ordinary differential equations (NODEs). NODEs have garnered immense interest for their desirable properties such as invertibility~\citep{Deng2020}, suitability for control~\citep{Asikis_2022}, and parameter efficiency and have been found more robust~\citep{Yan2019} than other deep neural network model variants.
Primarily used for developing black-box representations, NODEs can now be augmented to include domain priors and system structure via universal differential equations (UDEs)~\citep{Rackauckas2020}, making them an attractive proposition over other identification methods.
Hybrid neural networks form another class of flexible nonlinear models that combine a given, partially correct model with an additional neural network block to estimate the unmeasured, unmodeled dynamics that are difficult to capture from first principles or whitebox models \cite{psichogios1992hybrid}. This model is intuitive and easier to analyze, and has better prediction performance than the standard single blackbox or graybox models and has been applied in a number of domains \citep{meleshkova2021hybrid,gutierrez2011hybrid}.

In this work, we explore control-oriented NODE models for autonomous systems, such as AUVs, and study the effects of increasing physical system knowledge ranging from ``blackbox" models with no embedded physics knowledge to ``graybox" models which rely on varying degrees of physical knowledge. We consider six model types, including the vanilla NODE blackbox model, with and without domain priors and boundary constraints, the NODE graybox model, and multiple hybrid models with a fixed graybox block with varying parameter uncertainties augmented with a blackbox block that learns the residual errors. The key contributions include:
\begin{itemize}
    \item Development of a suite of domain-aware differentiable neural system models for AUVs
    \item Empirical analysis of the effect of applying domain-aware priors and system knowledge during training on open-loop prediction performance.
    
\end{itemize}
 We begin by introducing the ground truth AUV system model in Section \ref{sec:problemSetup}, which is used for data generation and model evaluation. Section \ref{sec: models} describes the different system models and the training procedure. We present results and discuss the model performance in Section \ref{sec:result}.

% System identification challenges (particularly in robotics) \\
% Physics-based modeling references in robotics \\
% Physics-based modeling limitations \\
% How ML/DL are used for modeling dynamical systems \\
% How they overcome the above limitations \\
% Add references of domain-aware modeling of dynamical systems in robotics\\
% Introduce Neural ODEs, advantages, citations \\

\section{Problem Setup}
\label{sec:problemSetup}

The problem addressed here is to develop a system identification approach for AUVs using NODEs that accommodates varying degrees of a priori information available for the system. Such system information considered here includes dynamical properties, such as state constraints and stability-related properties,  and system model structure. In this section, we present the dynamical model of an AUV as well as the development of the data set used for training in the system identification approach.
%\add{Jan: we need a short paragraph here to explain the modeling setup and assumptions on known model dynamics and data availability. We shall also mention challenges that we want to overcome with the used methods}

\subsection{AUV dynamic model}
We consider a torpedo-style vehicle, such as the Iver family of AUVs from \cite{iverWeb2022}, to generate data for our modeling and analysis. The nonlinear kinematic and hydrodynamic models we use are extended from  \cite{stankiewiczIver2021} and summarized here. 
\begin{subequations}
\begin{gather}
  \dot{p}_x = u\cos{\psi}\cos{\theta} + w\cos{\psi}\sin{\theta}, \\
  \dot{p}_y = u\sin{\psi}\cos{\theta} + w\sin{\psi}\sin{\theta}, \\
  \dot{p}_z = w\cos{\theta} - u\sin{\theta}, \\
  \dot{\theta} = q, \\
  \dot{\psi} = r/\cos{\theta}, \\
  \dot{u} = \overline{X}_{uu}u^2 + \overline{k}\delta_{u_c}, \\
  \dot{w} = \overline{Z}_{w|w|}w|w| + \overline{WB}\cos{\theta}, \\
  \dot{q} = \overline{M}_{uq}uq + \overline{M}_qq - \overline{Bz_B}\sin{\theta} + \bar{b}u^2\delta_{q_c}, \\
  \dot{r} = \overline{N}_{ur}ur + \bar{c}u^2\delta_{r_c} \\
  \dot{\delta}_{u_c} = \overline{K}_{\delta_u} (\delta_u - \delta_{u_c}) \\
  \dot{\delta}_{q_c} = \overline{K}_{\delta_q} (\delta_q - \delta_{q_c}) \\
  \dot{\delta}_{r_c} = \overline{K}_{\delta_r} (\delta_r - \delta_{r_c}).
\end{gather}
\label{eq:odes}
\end{subequations}
 The state is $\textbf{x}$ $=$ $(p_x, p_y, p_z, \theta, \psi, u, w, q, r,$ $\delta_{u_c},$ $\delta_{q_c},$ $\delta_{r_c})$, where $(p_x, p_y, p_z) \in \mathbb{R}^3$ are the respective north, east, down inertial coordinates, $\theta \in [-\pi/2, \pi/2], \psi \in \textsf{S}^1$ are the respective pitch and yaw orientation angles (Euler angles), $u, w \in \mathbb{R}$ are the respective surge and heave velocities, $ q, r \in \mathbb{R}$ are the pitch and yaw angular velocities, and $\delta_{u_c} \in [0, 1], \delta_{q_c}, \delta_{r_c} \in [-1, 1]$ are the delayed control inputs for thrust, elevator deflection, and rudder deflection, respectively. Note that $(u,w,q,r)\in \mathbb{R}^4$ are defined in the body frame. Similar to \cite{stankiewiczIver2021}, we assume roll and sway velocities are negligible. We also assume the vehicle has been trimmed to be neutrally buoyant and have no heave velocities, i.e. $w \equiv 0$. Refer to \cite{stankiewiczIver2021} for the hydrodynamic coefficients $(\overline{\cdot})$ of the Iver3 AUV. The control inputs $\textbf{u}=(\delta_u, \delta_q, \delta_r) \in \mathcal{U}:= [0,1] \times [-1,1]\times[-1,1]$ are the normalized thrust, deflections in the diving plane and rudder fins, respectively. The system \eqref{eq:odes} considers input delays that are common in practice as a low level controller tracks the given input commands. Here the input delay is modeled as a linear tracking control law as seen in the dynamics of $\delta_{u_c}, \delta_{q_c}$, and $\delta_{r_c}$ with $\overline{K}_{\delta_u} = \overline{K}_{\delta_q} =\overline{K}_{\delta_r} = -10$.
 
 A typical AUV hardware setup includes sensors like an IMU to measure angular rates, a compass for heading information, an acoustic doppler current profiler (ADCP) for measuring the vehicle's velocity through water and a depth sensor. We assume that the positions in the inertial coordinates are not directly measurable. Finally, the actuators contain sensors to determine the true actuator values in order to track the commanded input. Consistent with this realistic scenario, we consider the vehicle output to be velocities, angular rates and heading, $\mathbf{y}=(\theta, \psi, u, q, r, \delta_{u_c}, \delta_{q_c}, \delta_{r_c}) \in \mathcal{M}:=  [-\pi/2, \pi/2] \times  \mathbb{S}^1 \times \mathbb{R}^3 \times [0, 1] \times [-1, 1] \times [-1, 1]$, where $\mathbf{y}$ consists of $n_y = 8$ elements. 

% \in \mathbb{R}^3 \times \mathsf{S}^1 \times \mathsf{S}^1\times \mathsf{S}^1\times \mathbb{R}^6$,

\subsection{Dataset}
For the modeling process, we generate a dataset $\mathcal{D}$ of input-output pairs  obtained by applying appropriate excitation signals to explore the AUV state and input space:
\begin{equation}
    \mathcal{D}=\{[(\mathbf{u}_0, \mathbf{y}_0), (\mathbf{u}_{\delta}, \mathbf{y}_{\delta}), \ldots, (\mathbf{u}_{N\delta}, \mathbf{y}_{N\delta})]^{(i)}\},
\end{equation}
with sampling time $\delta \in \mathbb{R}$ and data batches $i=1, 2, \ldots, n$, of length $N$ each, and let $\mathbf{u}_k$ denote the vector $\mathbf{u}$ at time instant $k$. We use the mathematical model in \eqref{eq:odes} representing the AUV dynamics to generate the training and testing datasets and for ground truth trajectory generation for model evaluation. We sample control input trajectories and initial conditions from a distribution of pre-defined trajectories and initial conditions based on domain knowledge, that generates a combination of steady flight, diving, climbing, and port and starboard side turning maneuvers. Each input trajectory is constructed by combining a sequence of $50$ base trajectories for a length of $\frac{N}{50}$ seconds total per batch. Each base trajectory is randomly selected from a step response, periodic response, or spline constructed such that $\textbf{u}_{\tau}  \in \mathcal{U}$ for all $\tau \in [0, N\delta]$. The training set was determined by integrating \eqref{eq:odes} over a horizon of $N\delta$ s for $\delta = 0.01$ using `odeint' from Python's scipy module.

\section{Data-driven Models and Training}
\label{sec: models}

In this section we present the suite of models and how they are trained to approximate the true system \eqref{eq:odes}. The suite of models lie on a spectrum of domain-aware representations to reflect real-world scenarios depending on what degree of model knowledge is known a priori. The models are implemented as NODE.

\subsection{Blackbox model}

The blackbox model assumes minimal knowledge of the system dynamics. Here the dynamics are modelled with the neural network mapping $\textbf{f}_W: \mathcal{M} \times \mathcal{U} \to \mathbb{R}^{n_y}$, where $W$ represents the weights and biases of a multi-layer perceptron (MLP) with four hidden layers of size 128 each. The blackbox dynamical system is defined by:
\begin{equation}
    \dot{\textbf{z}}(t) = \textbf{f}_W(\textbf{z}(t), \textbf{u}(t)), \ \textbf{z}(0) = \textbf{z}_0 \in \mathcal{M},
\end{equation}
where $\mathbf{z} \in \mathbb{R}^{n_y}$ is the output of the model.

\subsection{Graybox model}

The graybox model requires the most knowledge of the system dynamics as it assumes the structure of \eqref{eq:odes} is known, but the model parameters are uncertain. Let the vector of model parameters be: $\boldsymbol{\mu} = (\overline{X}_uu, \overline{k},  \overline{M}_{uq}, \overline{M}_q, \overline{B_{z_B}}, \overline{b}, \overline{N}_{ur}, \overline{c} ) \in \mathbb{R}^{n_\mu}$, and let $\hat{\boldsymbol{\mu}} = (\hat{X}_uu, \hat{k}, \hat{M}_{uq}, \hat{M}_q, \hat{B_{z_B}}, \hat{b}, \hat{N}_{ur}, \hat{c})$ denote the approximate model parameters. Furthermore, let the approximate model dynamics be:
\begin{equation}\label{eq:graybox model}
    \textbf{f}_{\hat{\boldsymbol{\mu}}}(\textbf{y}, \textbf{u}) = \begin{pmatrix}
     q \\
   r/\cos{\theta} \\
  \hat{X}_{uu}u^2 + \hat{k}\delta_u \\
  \hat{M}_{uq}uq + \hat{M}_qq - \hat{Bz_B}\sin{\theta} + \hat{b}u^2\delta_q \\
  \hat{N}_{ur}ur + \hat{c}u^2\delta_r \\
  \dot{\delta}_{u_c} = \hat{K}_{\delta_u} (\delta_u - \delta_{u_c}) \\
  \dot{\delta}_{q_c} = \hat{K}_{\delta_q} (\delta_q - \delta_{q_c}) \\
  \dot{\delta}_{r_c} = \hat{K}_{\delta_r} (\delta_r - \delta_{r_c})
    \end{pmatrix}.
\end{equation}
The graybox dynamical system is defined as:
\begin{equation}
    \dot{\textbf{z}}(t) = \textbf{f}_{\hat{\boldsymbol{\mu}}}(\textbf{z}(t), \textbf{u}(t)), \ \textbf{z}(0) = \textbf{z}_0 \in \mathcal{M}.
\end{equation}

Assuming some degree of informed initial guess by the modeller, the graybox model parameters, $\hat{\boldsymbol{\mu}}$, are initialized to within an offset of the true model parameters, $\boldsymbol{\mu}$. Let $e_\mu \in [0, 1]$ represent the error in the model parameters. For our experiments, parameters are sampled uniformly from a set whose range is defined by $e_\mu$. For example, $\hat{X}_{uu} \sim \mathscr{U}([\min \{X_{uu}(1-e_\mu), X_{uu}(1+e_\mu)\}, \max \{X_{uu}(1-e_\mu), X_{uu}(1+e_\mu) \}]$, where $\mathscr{U}(S)$ denotes a uniform sampling of the set $S$. To state this more formally, let $\mathcal{S}^\mu_{e_\mu} = \mathcal{S}^{\mu^1}_{e_\mu} \times ...\times \mathcal{S}^{\mu^{n_\mu}}_{e_\mu}$, where $\mathcal{S}^{\mu^i}_{e_\mu} = [\min \{\mu_i(1-e_\mu), \mu_i(1+e_\mu)\}, \max \{\mu_i(1-e_\mu),\mu_i(1+e_\mu) \} ]$ with $\mu_i$ as the $i$th element of $\boldsymbol{\mu}$. Thus for a given ${e_\mu}$, $\hat{\boldsymbol{\mu}} \sim \mathscr{U}(\mathcal{S}^\mu_{e_\mu})$ implies an estimated parameter uniformly sampled from a range $e_\mu$ of the true model parameters. For the graybox model, the model parameters were initialized as $\hat{\boldsymbol{\mu}}_0 \sim \mathcal{S}^\mu_{e_\mu=1.0}$ to provide the maximum offset of the true values and test if the estimated model parameters can converge to the true model parameters.

\subsection{Hybrid model}

The hybrid model combines an uncertain approximation of the system dynamics, similar to the graybox, with a blackbox model to  compensate for residual errors. This model considers the case when some, imperfect knowledge, of the system is known a priori. Here the model parameters of the graybox model, $\hat{\boldsymbol{\mu}}$, are initialized and fixed to a random offset of their true parameters. The approximate model used is the same as the graybox $\textbf{f}_{\hat{\boldsymbol{\mu}}}(\textbf{y}, \textbf{u})$ from \eqref{eq:graybox model}. The blackbox model is the same as from the previous case with the neural network mapping $\textbf{f}_W: \mathcal{M} \times \mathcal{U}\to \mathbb{R}^{n_y}$, with $W$ representing the weights and biases in the neural network. In this hybrid case, only the blackbox parameters, $W$, are trained and not the offset model parameters $\hat{\boldsymbol{\mu}}$. The hybrid model dynamics are defined by:
\begin{equation}
    \dot{\textbf{z}}(t) = \textbf{f}_{\hat{\boldsymbol{\mu}}}(\textbf{z}(t), \textbf{u}(t)) + \textbf{f}_W(\textbf{z}(t)), \textbf{u}(t)), \ \textbf{z}(0) = \textbf{z}_0 \in \mathcal{M}.
\end{equation}

To compare various different errors in the graybox portion of the hybrid model, we look at three implementations of the hybrid model which are dependent on the percentage of error in the model parameters, $\hat{\boldsymbol{\mu}}$, which is defined by $e_\mu$. Here we look at the three cases `$\text{hybrid}_{e_\mu = 1.0}$', `$\text{hybrid}_{e_\mu = 0.5}$', and `$\text{hybrid}_{e_\mu = 0.3}$'.

\subsection{Training}
We perform batched training on the trajectories in $\mathcal{D}$. Each model was trained for $n = 5$ batches. The length, $N$, was increased by a factor of 2 for each batch yielding batches with $N=100$, $N = 200$, $N = 400$, $N = 800$, and $N = 1600$. The purpose of using multiple batches of different lengths is to prevent over-approximation over the same time length and it extends how long the control inputs are applied to the system, which yields different system responses. Training is performed using the Neuromancer differentiable programming framework \cite{Neuromancer2022}, wherein the following optimization problem is solved to train the model parameters $\bf V \in\mathbb{R}^{n_v}$:
% {\bf z}_{k+1}^i = {\bf z}_{k}^i + \delta {\bf f}_{\bf V^i}( {\bf z}_k^i,{\bf u}_k),
\begin{subequations}
\label{eq:neuromancer}
    \begin{align}
 \min_{{\bf V}} & \frac{1}{mN} \sum_{i=1}^{m} \sum_{k=0}^{N-1}  \big( \ell_{\texttt{mse}}( {\bf z}_k^i, {\bf y}_k) +   p_z(c({\bf z}_k^i)   \big) & 
 \label{eq:DPC:objective} \\ 
  \text{s.t.} \ &  {\bf z}_{k+1}^i = ODESolve({\bf f}_{\bf V^i}( {\bf z}_k^i,{\bf u}_k)), \  k \in \mathbb{N}_{0}^{N-1}   \label{eq:dpc:x}  \\
   \ &   {\bf z}_{0}^i = {\bf y}_0 \label{eq:dpc:initial}  \\
  \ & {\bf V}^i  \in \Xi \subset \mathbb{R}^{n_v} \label{eq:dpc:xi}
\end{align}
\end{subequations}
where 
\begin{equation}\label{eq:loss objective}
    \ell_{\texttt{mse}}( {\bf z}, {\bf y}) = \| {\bf y} -  {\bf z} \|_2^2
\end{equation}
is the objective loss function, $c:\mathbb{R}^{n_y} \to \mathbb{R}^{n_c}$ represents state constraints 
with the penalty functions  $p_x: \mathbb{R}^{n_c} \to \mathbb{R} $, $\Xi$ is a random distribution from which the model parameters ${\bf V}^i$ are sampled from, and $\textbf{f}_{\textbf{V}}: \mathcal{M} \times \mathcal{U} \to \mathbb{R}^{n_y}$ is a general NODE model parameterized by $\textbf{V}$. Stochastic gradient descent is used to minimize the loss function \eqref{eq:DPC:objective}
over the distribution of model parameters \eqref{eq:dpc:xi}, where $m$ represents the total number of model parameter samples, and $i$ denotes the index of the sample. We note that ODESolve can represent any ODE solver such as Runge Kutta or Euler integrators. In our implementation, an Euler approximation is used to integrate the NODE models. The problem \eqref{eq:neuromancer} is solved for each batch, and the final trained model is output after the final batch. The penalty functions are used to enforce boundary conditions, which will be discussed in the next section.

We emphasize that the different model types will have different model parameters that are trained. In the case of the blackbox model, $\textbf{f}_{\textbf{V}} = \textbf{f}_{\textbf{W}}$, where $\textbf{V} = \textbf{W}$ are the MLP parameters that will be trained. For the graybox model, $\textbf{f}_{\textbf{V}} = \textbf{f}_{\hat{\boldsymbol{\mu}}}$, where $\textbf{V} = \hat{\boldsymbol{\mu}}$. Finally, for the hybrid model $\textbf{f}_{\textbf{V}} = \{\textbf{f}_{\hat{\boldsymbol{\mu}}}, \textbf{f}_W$\}, where $\textbf{V} = \textbf{W}$.

\subsection{Learning with Boundary Conditions}

The Neuromancer framework allows for enforcing boundary conditions in the form of state constraints $c({\bf z}) $. In the case of \eqref{eq:odes}, it is clear that if $\theta \to \pm \frac{\pi}{2}$, then the system dynamics become undefined at the singularity. We note that the randomized training sets are restricted to $\mathcal{M}$ so such instabilities are never reached, and thus we aim to enforce constraints on the models to also avoid the singularity. Similarly, the states associated with input delays should remain inside the input set $\mathcal{U}$. 

To enforce the boundary conditions, we define $c_i(\mathbf{z}) := \max \{0, -\mathbf{z}_i + \underline{z}_i \} + \max \{0, \mathbf{z}_i - \overline{z}_i \}$, where $i \in \mathbb{N}_0^{n_y -1}$, $\mathbf{z}_i$ is the $i$th element of $\mathbf{z}$, $\underline{z}$ $=$ $( -\frac{\pi}{2},\linebreak[1] -\infty,\linebreak[1] -\infty,\linebreak[1] -\infty,\linebreak[1] -\infty,\linebreak[1] 0,\linebreak[1] -1,\linebreak[1] -1)$ is the lower bound on $\mathbf{z}$, and $\overline{z} = (\frac{\pi}{2}, \infty,\infty,\infty,\infty,1, 1, 1)$ is the upper bound on $\mathbf{z}$. Furthermore, let $c(\mathbf{z}) = (c_0(\mathbf{z}), ..., c_{n_y-1}(\mathbf{z}))$. The penalty function used is: $p_c(\textbf{c}) = \| \textbf{c} \|_2$.

\subsection{Learning with Domain-aware Priors}

In addition to state constraints, the Neuromancer framework enables the training process to enforce system-type properties of the trained model. As with many physical systems, the AUV system is dissipative, which means that the system has no internal production of energy, and the energy stored by the system at any time does not exceed the energy externally supplied at that time \citep{willems2007dissipative}. The trained models can be trained to enforce the this property. For the blackbox model and the blackbox component of the hybrid model, if the 2-norm of the neural network weights and activation matrices are restricted to be less than $1$, then the NODE is dissipative \citep{Drgona2022a}. Here, the neural network weights are restricted such that the eigenvalues of the weighting and activation matrices lie within $[\sigma_{min}, \sigma_{max}]$ for $\sigma_{min}, \sigma_{max} \in \mathbb{R}_{\geq 0}$. This restriction is applied during training to always ensure the dissipativity property.

For the blackbox model, the eigenvalues are set to $\sigma_{min} = 0.5$, $\sigma_{max} = 1.0$. The reason for setting $\sigma_{min} = 0.5$ instead of zero is to keep the model away from steady-state solutions since the training set was designed to excite the system and avoid such solutions. The blackbox component of the hybrid models had eigenvalues set to $\sigma_{min} = 0.0$, $\sigma_{max} = 1.0$. Here we allow the eigenvalues to reach zero unlike in the blackbox model because in the hybrid model the blackbox component is used to compensate for the residual error of the graybox. Thus if the graybox is able to significantly reduce the residual error, the blackbox should approximate a zero map. To show the effects of the boundary conditions and the domain-aware priors on the neural network NODE model, we compare in our experiments a NODE model without boundary conditions/domain-aware priors denoted ``blackbox", and the same architecture with boundary conditions/domain-aware priors denoted ``constr. blackbox"

%\add{Jan:  we can have domain structural priors as separate section describing how to modify all our baseline black box models NODEs, NSSMs, and Koopmans to include priors of the IVER vehicle type models \\
%gray box parameter estimation in NODE \\
%eigenvalue, sparsity constraints in NSSM, Koopman }

% \subsection{Block SSM}

% \begin{figure}[h]
%     \centering
%     \includegraphics[scale=0.6]{Figures/blockSSM_best.png}
%     \caption{Best BlockSSM, best loop dev Y pred dynamics eq Yf  = 0.003 }
%     \label{fig:my_label}
% \end{figure}

%\begin{figure}[h]
%    \centering
%    \includegraphics[scale=0.6]{Figures/neural_ode_blackbox_best.png}
 %   \caption{Best Neural ODE Blackbox, eval %mse  = 0.0009757 }
%    \label{fig:my_label}
%\end{figure}

%===============================================================================

\section{Results and Discussion}
\label{sec:result}

\begin{figure*}
\begin{center}
    \centering
    \subfigure[Pitch angle ($\theta$) residual MSE]{
        \centering
        \includegraphics[width=5.18cm]{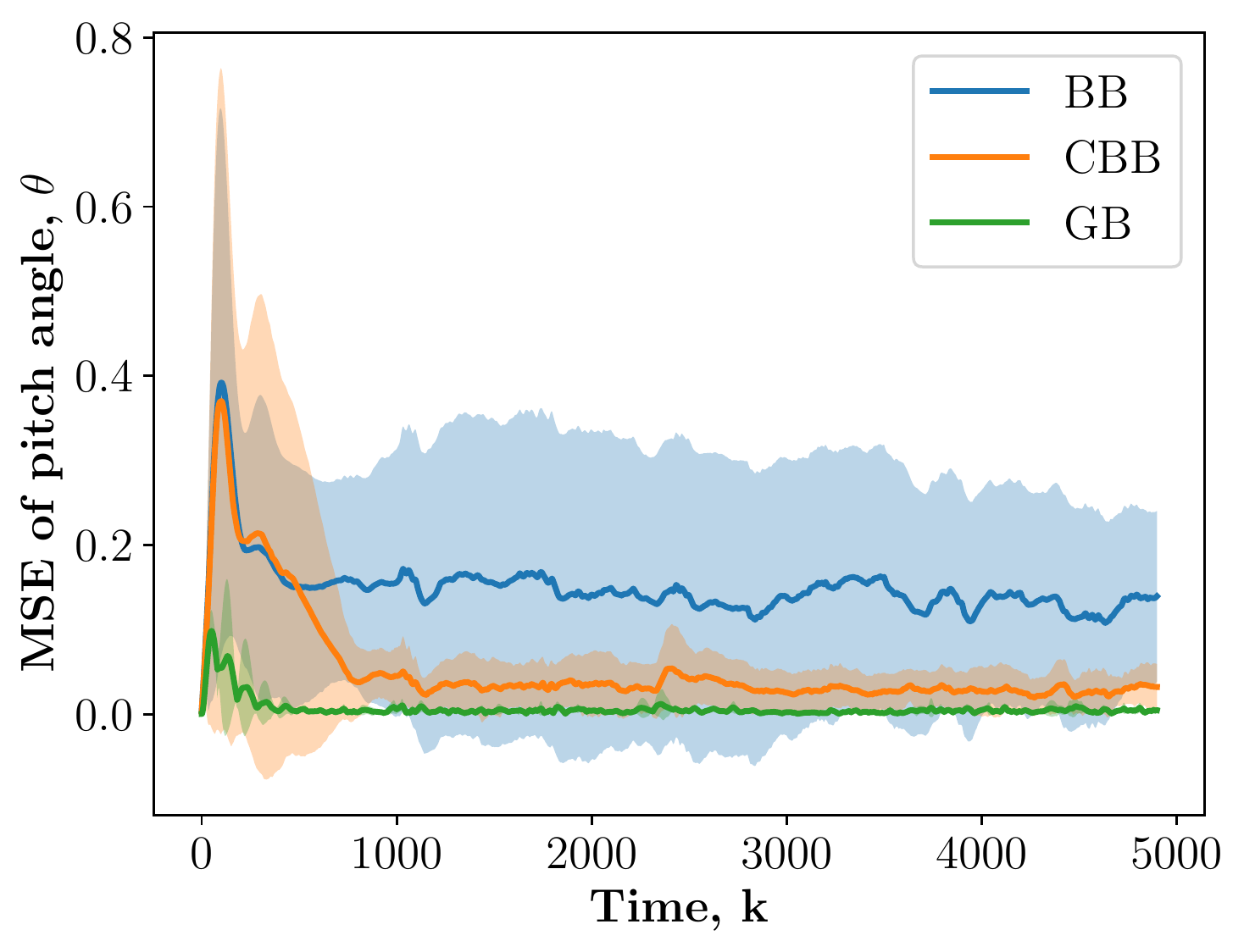}
        \includegraphics[width=5.22cm]{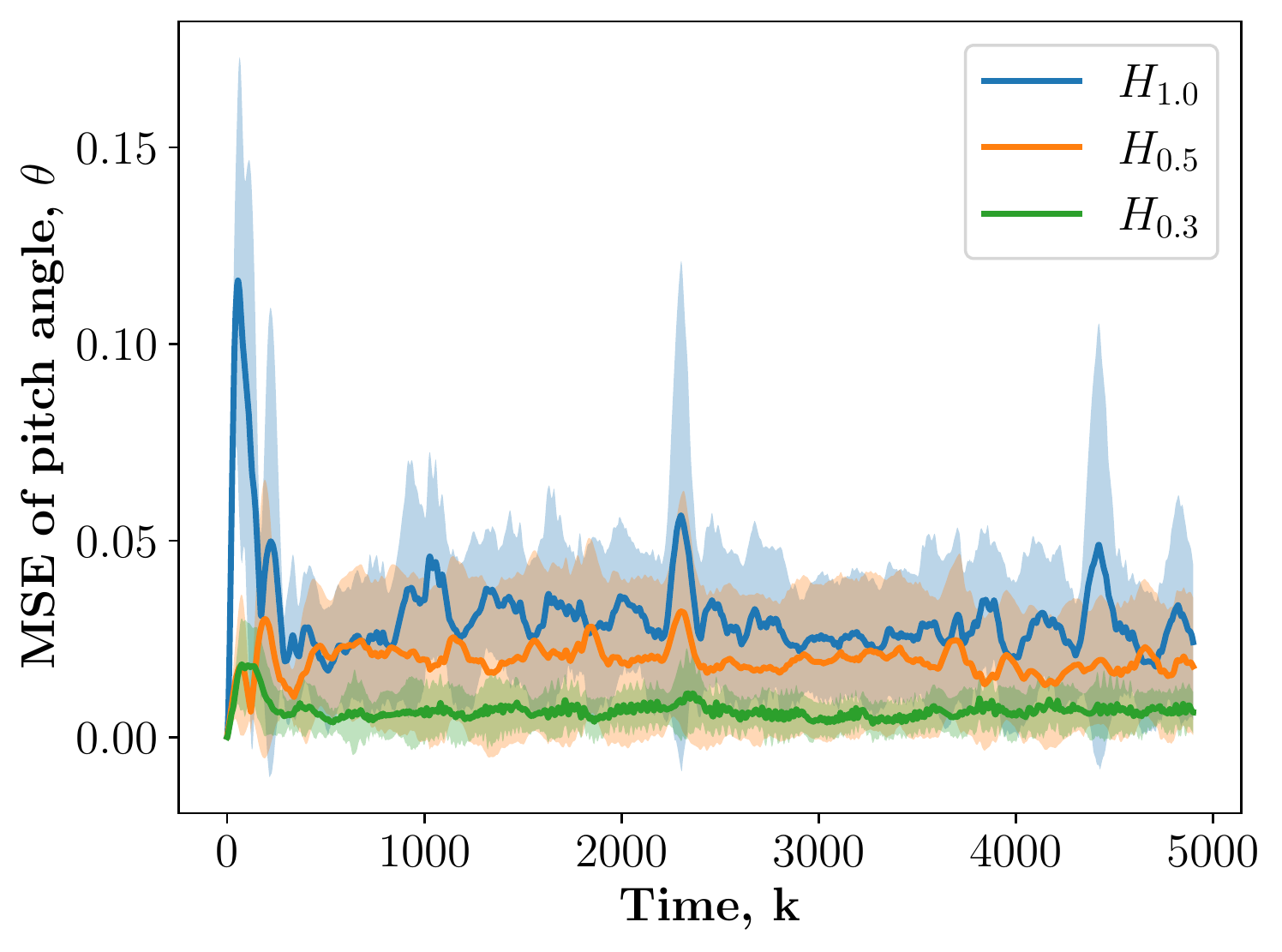}
    }
        \subfigure[constr. blackbox (CBB) ]{
        \centering
        \includegraphics[width=5.7cm]{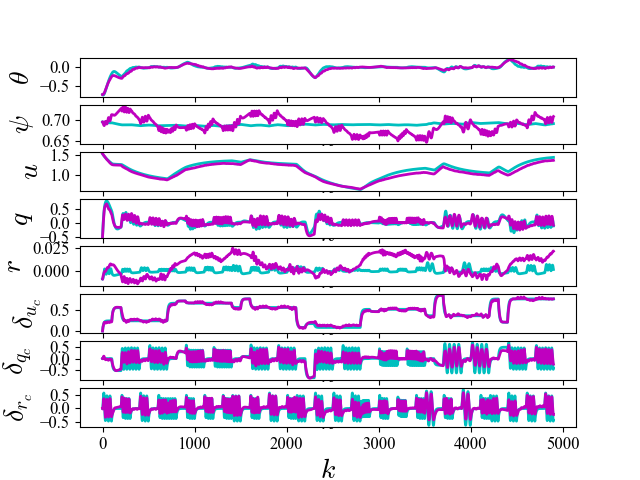}
    }
    \subfigure[Surge velocity ($u$) residual MSE]{
        \centering
        \includegraphics[width=5.2cm]{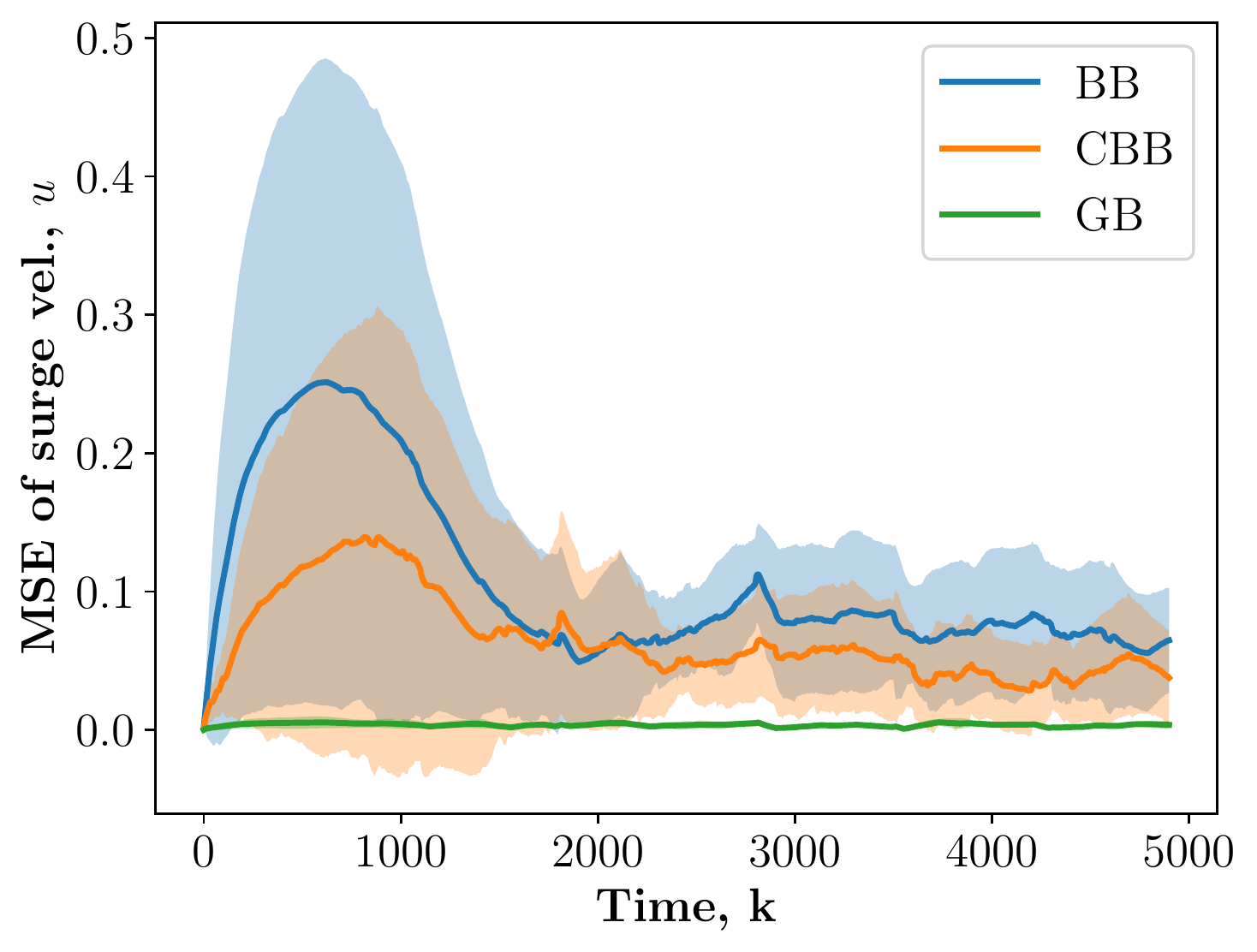}
        \includegraphics[width=5.2cm]{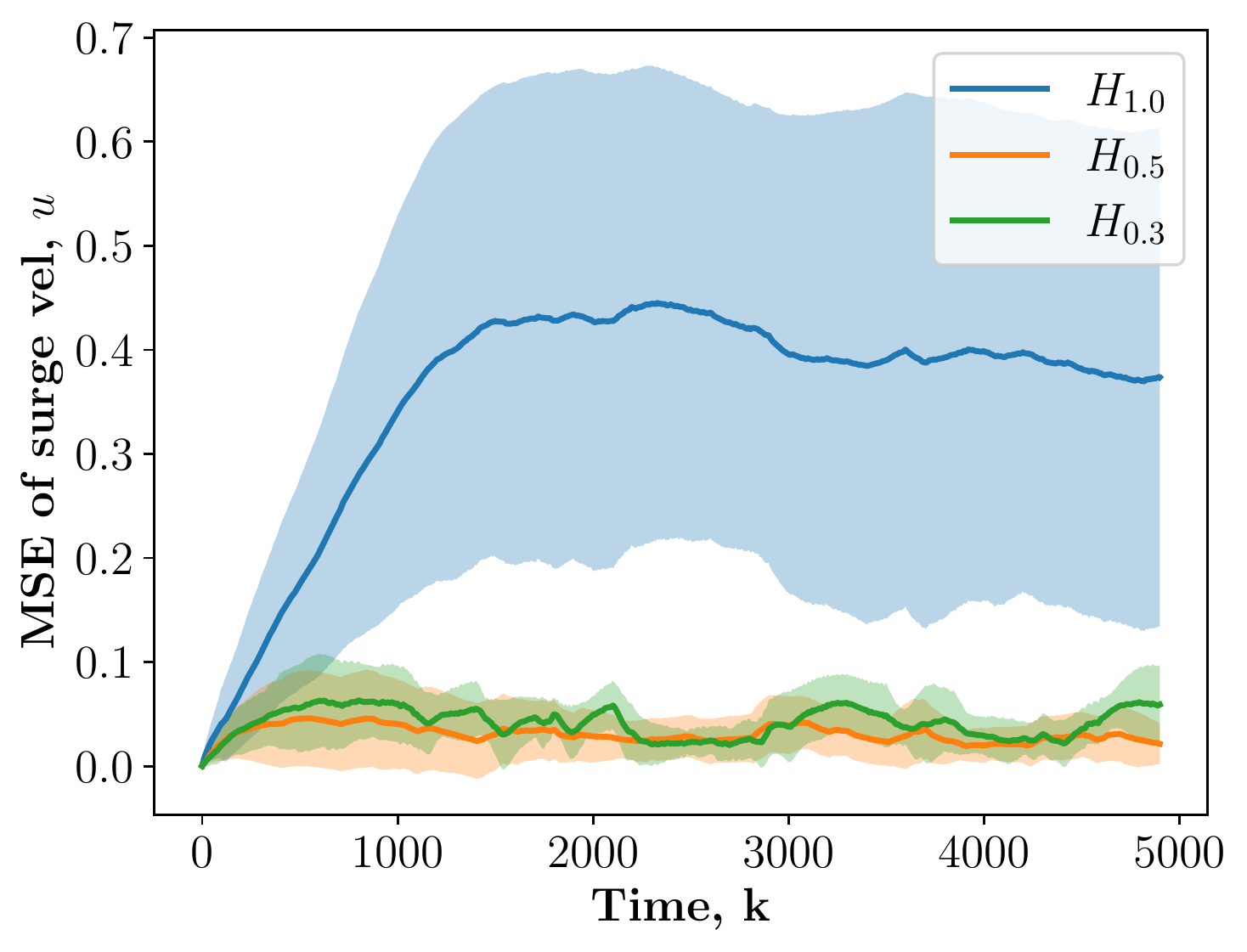}
    }
    \subfigure[graybox (GB)]{
    \centering
    \includegraphics[width=5.7cm]{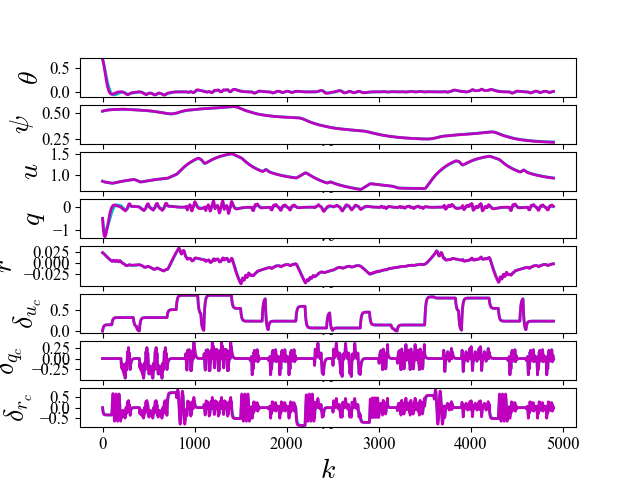}
    }

    \caption{Comparison of MSE for (a) Pitch angle, $\theta$, (c) Surge velocity, $u$, for each model type. The open-loop response, showing the measured states of the best model of type (b) constrained blackbox (CBB) and (d) graybox (GB), evaluated on test data. The ground truth and the predicted trajectories are in cyan and purple, respectively. }
    \label{fig:sim_results}
\end{center}
\end{figure*}

For our empirical analysis, 10 sets of randomly sampled initial parameters for each of the six model types, blackbox, constr. blackbox , $\text{hybrid}_{e_\mu = 1.0}$, $\text{hybrid}_{e_\mu = 0.5}$, $\text{hybrid}_{e_\mu = 0.3}$, and graybox, were fit to the training dataset using the ADAM-W~\citep{loshchilov2017decoupled} variant of stochastic gradient descent. 
% We note that the stochastic nature of the Neuromancer framework yields different results for the same initial condition. 
% To account for variations in the trained models, each model type was trained 10 times for the same initial conditions, thus yielding 10 models per model type. 
After training these 10 models instances per model type, they were evaluated on a held out test set. The test set was developed by randomly initializing five initial conditions from $\mathcal{M}$ each for five randomly generated input trajectories, which yields a total of 25 trajectories. The input trajectories were developed in the same manner as the training set with $N = 5000$. We note that the test set is the same for all model types.

For each model type, there are 10 models that were each evaluated on 25 trajectories. For each trajectory of each model, the mean squared error (MSE) is:
\begin{equation}
    \frac{1}{N} \sum_{k = 0}^{k = N} \ell_{\texttt{mse}}( {\bf z}'_k, {\bf y}'_k)
\end{equation}
where ${\bf z}'_k, {\bf y}'_k$ are the normalized outputs of the model and true system respectively. The outputs are normalized using the training set data for a unit-less error metric. 
% The mean value over all 10 models of the mean of the mse across all 25 trajectories is referred to as the `mean of all'. The associated standard deviation of the `mean of all', is referred to as the `std of all'.
We note that the inter quartile range approach was used to remove outlier models from the results. For each model type, we consider the `best' trained model to be the one in which the mean of the mean square error (MSE) over all 25 trajectories is the lowest.
% The mean of the mse of the best model over all 25 trajectories is referred to as the `mean of best', and the associated standard deviation is the `std of best'. Finally, the `best' trajectory of the best model is one of the 25 trajectories where in the best model yields the lowest mse. The mse associated with the best trajectory of the best model is referred to as the `best mse.' 

\begin{table}[hbt!]
%\begin{center}
\centering
\caption{Comparison of blackbox plain, blackbox, hybrid, and graybox models  }
\label{tab:model comparison}
\begin{tabular}{llllll}
\toprule
  \textbf{model} &  \textbf{Mean MSE} &  \textbf{Standard deviation} \\
\midrule
blackbox &    16.459841 &   61.742187 \\
constr. blackbox &      7.985731 &    0.311266 \\
   $\text{hybrid}_{e_\mu = 1.0}$  &  82.699305 & 11399.445734  \\
   $\text{hybrid}_{e_\mu = 0.5}$  &   6.01626 &   21.648181 \\
   $\text{hybrid}_{e_\mu = 0.3}$  &    2.626124 &    0.474631\\
   graybox &     0.025351 &    0.000446 \\
\bottomrule
\end{tabular}
%\end{center}
\end{table}

Table \ref{tab:model comparison} shows the results of the blackbox, constr. blackbox , $\text{hybrid}_{e_\mu = 1.0}$, $\text{hybrid}_{e_\mu = 0.5}$, $\text{hybrid}_{e_\mu = 0.3}$, and graybox models. 
Figure \ref{fig:sim_results} interprets Table \ref{tab:model comparison} visually and focuses on the pitch and surge velocity residual error evolution over time and the standard deviation over the 25 trajectories, for each learned model type. The open loop response of all measured states for the constr. blackbox and the graybox models are also presented.
% We note that the evaluation trajectory associated with the best MSE is not consistent across all model types. 

Several conclusions can be drawn from Table \ref{tab:model comparison}.  First, the graybox outperforms all model types with respect to the mean of all trained models and associated standard deviation. This shows that on average, the graybox is the best model type to use and the trained models will consistently yield the same performance in terms of mean squared error from the true model dynamics. This fits with intuition in that the graybox requires the most domain information to describe the system structure exactly, but with uncertain parameters. However this may not be true in many scenarios especially when simplifications are made in the graybox model. 

In practice, there usually will be unmodeled dynamics that need to be accounted for. This can be represented by the hybrid model, where imperfect model parameters will cause the graybox system dynamics to deviate with the true system and requires the blackbox component to compensate. The hybrid model results from Table \ref{tab:model comparison} also follow intuition in that as $e_\mu \to 0.0$, i.e., the graybox component approaches the true system, the blackbox component performs better with respect to reducing the MSE between the hybrid model and the true system model. Furthermore, the trained models become more consistent in their performance as $e_\mu \to 0.0$, which is represented by the decrease in standard deviation of the model performances.

The blackbox/blackbox constrained models represent the case when no model knowledge is assumed and they must try to learn the true system dynamics. According to Table \ref{tab:model comparison}, the blackbox and blackbox constrained models actually outperform the $\text{hybrid}_{e_\mu = 1.0}$ in terms of mean MSE and standard deviation of all trained models. This suggests that the blackbox type models are preferred to the hybrid model if the graybox component poorly represents the true system. In this experiment, a poor representation is defined by $e_\mu = 1.0$, wherein the model parameters are sufficiently inaccurate as to impede the blackbox component's ability to compensate for the error between the graybox and true system dynamics. One possible reason for this is that the model parameters may cause the graybox model to approach an unstable system dynamics, which would mean that the blackbox component is attempting to compensate for the near-instability while also eliminating the residual error between the graybox and true system.

Finally, we can see the comparison between blackbox and blackbox constrained models. It is clear that on average, the blackbox constrained models yield lower MSE and are more consistent with respect to the standard deviation of all models. This suggests that the use of a priori information of the system model, specifically the state constraints and dissipativity of the system, help improve the training of the blackbox models and can reduce the errors of the blackbox model by a factor of 2. This result fits with intuition in that more information leads to improved system identification.

%===============================================================================

\section{Conclusion}
\label{sec:conclusion}
In this paper we present a domain-aware framework for system identification of autonomous underwater vehicles using neural ordinary differential equations (NODEs). The framework allows for structured training of the NODE model by enforcing state constraints and domain-aware priors to ensure dissipativity of the closed-loop system. The proposed approach considers a suite of models that lie on a spectrum of domain-aware representations to reflect real-world scenarios depending on what degree of model knowledge is known a priori. The models range from a blackbox model,  a hybrid model, and a graybox model. The blackbox, graybox, and several implementations of the hybrid model were compared in an experiment using training sets and evaluation sets.
%The results show that the graybox outperforms all models, as expected, but the blackbox may outperform the hybrid model if the a priori system knowledge is sufficiently incorrect. 
Future work will investigate the effect of noise on training and incorporate the domain-aware models in a predictive control structure.

\bibliography{references}  

\end{document}